\documentclass[12pt,a4paper,twoside]{article}
\usepackage{a4wide,amsmath,amsbsy,amsfonts,amssymb,stmaryrd,amsthm}
\newcommand\ZZ{{\hat{\mathbb Z}}}
\newcommand\Z{{\mathbb Z}}
\newcommand\PP{{\mathbb P}}

\newcommand\C{{\mathbb C}}

\newcommand\ra{\rightarrow}

\newcommand\SL{\operatorname{SL}}

\newcommand\Sp{\operatorname{Sp}}
\newcommand\aut{\operatorname{Aut}}
\newcommand\out{\operatorname{Out}}

\newcommand\diff{\operatorname{Diff}}
\newcommand\lra{\longrightarrow}
\newcommand\hookra{\hookrightarrow}
\newcommand\tura{\twoheadrightarrow}
\newcommand\da{\downarrow}

\renewcommand{\hom}{\operatorname{Hom}}

\newcommand\sr{\stackrel}
\newcommand\st{\scriptstyle}
\newcommand\sst{\scriptscriptstyle}
\newcommand\cGG{\check{\GG}}
\newcommand\hGG{\hat{\GG}}
\newcommand\hH{\hat{H}}

\newcommand\hP{\hat{\Pi}}

\newcommand\tGG{\tilde{\GG}}

\newcommand\ssm{\smallsetminus}
\newcommand\ol{\overline}

\newcommand\cM{{\cal M}}
\newcommand\cH{{\cal H}}

\newcommand\cC{{\cal C}}

\newcommand\GG{\Gamma}
\newcommand\ld{\lambda}

\newcommand\gm{\gamma}

\def\co{\colon\thinspace}

\newtheorem{theorem}{Theorem}[section]
\newtheorem{corollary}[theorem]{Corollary}
\newtheorem{proposition}[theorem]{Proposition}
\newtheorem{lemma}[theorem]{Lemma}

\theoremstyle{definition}
\newtheorem{definition}[theorem]{Definition}   
\newtheorem{remark}[theorem]{Remark}

\begin{document}

\title{The congruence subgroup property for the \\ hyperelliptic modular group: the open surface case}
\author{Marco Boggi}\maketitle

\begin{abstract}Let ${\cal M}_{g,n}$ and ${\cal H}_{g,n}$, for $2g-2+n>0$, be, respectively, the 
moduli stack of $n$-pointed, genus $g$ smooth curves and its closed substack consisting
of hyperelliptic curves. Their topological fundamental groups can be identified, respectively, 
with $\Gamma_{g,n}$ and $H_{g,n}$, the so called {\it Teichm{\"u}ller modular group} and 
{\it hyperelliptic modular group}. A choice of base point on ${\cal H}_{g,n}$ defines a 
monomorphism $H_{g,n}\hookrightarrow\Gamma_{g,n}$. 

Let $S_{g,n}$ be a compact Riemann surface of genus $g$ with $n$ points removed.
The Teichm\"uller group $\Gamma_{g,n}$ is the group of isotopy classes of diffeomorphisms 
of the surface $S_{g,n}$ which preserve the orientation and a given order of the punctures.
As a subgroup of $\Gamma_{g,n}$, the hyperelliptic modular group then admits a natural faithful 
representation $H_{g,n}\hookrightarrow\operatorname{Out}(\pi_1(S_{g,n}))$.

The {\it congruence subgroup problem for} $H_{g,n}$ 
asks whether, for any given finite index subgroup $H^\lambda$ of $H_{g,n}$, 
there exists a finite index characteristic subgroup $K$ of $\pi_1(S_{g,n})$
such that the kernel of the induced representation $H_{g,n}\to\operatorname{Out}(\pi_1(S_{g,n})/K)$ 
is contained in $H^\lambda$. The main result of the paper is an
affirmative answer to this question for $n\geq 1$.
\newline


\end{abstract}    


\section{Introduction}\label{intro}
Let $S_{g,n}$, for $2g-2+n>0$, be the differentiable surface obtained from a compact Riemann surface
$S_g$ of genus $g$ removing $n$ distinct points $P_i\in S_g$, for $i=1,\ldots,n$. The Teichm\"uller
modular group of $S_{g,n}$ is defined to be the group of isotopy classes of diffeomorphisms or, 
equivalently, of homeomorphisms of the surface $S_{g,n}$ which preserve the orientation and the
given order of the punctures:
$$\GG_{g,n}:=\diff^+(S_{g,n})/\diff_0(S_{g,n})\cong\hom^+(S_{g,n})/\hom_0(S_{g,n}),$$
where $\diff_0(S_{g,n})$ and $\hom_0(S_{g,n})$ denote the connected components of the identity 
in the topological groups of diffeomorphisms $\diff^+(S_{g,n})$ and of
homeomorphisms $\hom^+(S_{g,n})$.

Let $\Pi_{g,n}$ denote the fundamental group of $S_{g,n}$ for some choice of base point. 
From the above definition and some elementary topology, it follows that there is a faithful 
representation:
$$\rho\co \GG_{g,n}\hookra\out(\Pi_{g,n}).$$
A {\it level} of $\GG_{g,n}$ is just a finite index subgroup $H<\GG_{g,n}$.
A characteristic finite index subgroup $\Pi^\ld$ of $\Pi_{g,n}$ determines the {\it geometric level}
$\GG^\ld$, defined to be the kernel of the induced representation:
$$\rho_\ld\co \GG_{g,n}\lra\out(\Pi_{g,n}/\Pi^\ld).$$
The {\it congruence subgroup problem} asks whether geometric levels are cofinal in the set
of all finite index subgroups of $\GG_{g,n}$, ordered by inclusion.

This problem is better formulated in the geometric context of moduli spaces of curves.
Let $\cM_{g,n}$, for $2g-2+n>0$, be the moduli stack of $n$-pointed, genus $g$, smooth algebraic 
complex curves. It is a smooth connected Deligne-Mumford stack (briefly {\it D-M stack}) over $\C$ of
dimension $3g-3+n$, whose associated underlying complex 
analytic and topological {\'e}tale groupoids, we both denote by $\cM_{g,n}$ as well.

In the category of analytic {\'e}tale groupoids, there are natural
and general definitions of topological homotopy groups (see \cite{Noohi}). However, 
for stacks of the kind of $\cM_{g,n}$, such groups can be described in a 
simpler way. In fact, $\cM_{g,n}$ has a universal cover $T_{g,n}$ in the category of analytic
manifolds. The fundamental group $\pi_1(\cM_{g,n},[C])$ is then identified with the deck 
transformations' group of the cover $T_{g,n}\ra\cM_{g,n}$ and the higher homotopy groups are 
naturally isomorphic to those of $T_{g,n}$.

From this perspective, Teichm{\"u}ller theory is the study of the geometry of the universal cover 
$T_{g,n}$ of the moduli space $\cM_{g,n}$, called Teichm{\"u}ller space, and of its topological 
fundamental group $\pi_1(\cM_{g,n},[C])$. The basic facts of Teichm\"uller theory are that 
$T_{g,n}$ is contractible, thus making of $\cM_{g,n}$ a classifying space for $\GG_{g,n}$, and
that the choice of a lift of a point $[C]\in\cM_{g,n}$ to $T_{g,n}$ and of a diffeomorphism 
$S_{g,n}\ra C\ssm\{\mbox{marked points}\}$ identifies the Teichm\"uller modular group $\GG_{g,n}$ 
with $\pi_1(\cM_{g,n},[C])$. The representation:
$$\rho\co\pi_1(\cM_{g,n},[C])\ra\out(\Pi_{g,n}),$$
induced by the identification of $\GG_{g,n}$ with $\pi_1(\cM_{g,n},[C])$,
is equivalent to the universal topological monodromy representation associated with the universal punctured curve $\cM_{g,n+1}\ra\cM_{g,n}$. Algebraically, this may be recovered as the outer 
representation associated to the short exact sequence determined on topological fundamental groups 
by this curve:
$$1\ra\Pi_{g,n}\ra\pi_1(\cM_{g,n+1})\ra\pi_1(\cM_{g,n})\ra 1.$$

The algebraic fundamental group of a D-M stack $X$ over $\C$ is naturally isomorphic to the 
profinite completion $\hat{\pi}_1(X)$ of its topological fundamental group $\pi_1(X)$. 
It basically follows from the triviality of the center of the profinite completion $\hP_{g,n}$ of $\Pi_{g,n}$ 
that the above fibration induces on algebraic 
fundamental groups the short exact sequence:
$$1\ra\hP_{g,n}\ra\hat{\pi}_1(\cM_{g,n+1})\ra\hat{\pi}_1(\cM_{g,n})\ra 1.$$
The {\it universal algebraic monodromy representation} is the outer representation:
$$\hat{\rho}\co\hat{\pi}_1(\cM_{g,n})\ra\out(\hP_{g,n}),$$
associated to the above short exact sequence. It is not hard to see that the congruence subgroup property holds for $\GG_{g,n}$ if and only if the representation $\hat{\rho}$ is faithful.

In \cite{B1}, a positive answer to the above question was claimed but a gap emerged in an essential
step of the proof (more precisely, in the proof  of Theorem~5.4). As it is explained in 
detail below, this paper recovers some of the results of \cite{B1}. 

Indeed, the congruence subgroup problem can be formulated for any special subgroup of
the Teichm\"uller group. The case we will deal with in this paper is that of the fundamental 
group of the closed sub-stack $\cH_{g,n}$ of
$\cM_{g,n}$ parametrizing smooth hyperelliptic complex curves, for $g\geq 1$. 
Observe that, for $g=1, 2$, all curves are hyperelliptic, i.e. admit a degree $2$ morphism onto
$\PP^1$. We then define the {\it hyperelliptic modular group} to be the topological 
fundamental group of the stack $\cH_{g,n}$.

It is a classical fact of Teichm\"uller theory that the subspace of the Teichm\"uller space $T_{g,n}$,
parametrizing hyperelliptic curves, consists of a disjoint union of contractible analytic subspaces.
The natural embedding $\cH_{g,n}\subset\cM_{g,n}$ then induces, choosing for base points the 
isomorphism class $[C]$ of a hyperelliptic curve, a monomorphism of topological fundamental groups 
$\pi_1(\cH_{g,n},[C])\hookra\pi_1(\cM_{g,n},[C])$. 
Let us remark that the image of the latter map, in general, is not a normal subgroup of 
$\pi_1(\cM_{g,n},[C])$. 

After the identification of $\pi_1(\cM_{g,n},[C])$ with $\GG_{g,n}$, we denote the subgroup 
corresponding to $\pi_1(\cH_{g,n},[C])$ simply by $H_{g,n}$. Let then $\iota$ be the element 
of $\GG_{g,n}$ corresponding to the hyperelliptic involution on $C$.
For $g\geq 2$ and $n=0$ or $g=1$ and $n=1$, the subgroup $H_{g,n}$ is 
the centralizer of $\iota$ in $\GG_{g,n}$. 

For a given characteristic subgroup of finite index $\Pi^\ld$ of $\Pi_{g,n}$, let us
define $H^\ld:=H_{g,n}\cap\GG^\ld$ and call it the geometric level of $H_{g,n}$ associated to $\Pi^\ld$. 
The {\it congruence subgroup problem for the hyperelliptic modular group} asks whether
geometric levels of $H_{g,n}$ are cofinal in the set of finite index subgroups of $H_{g,n}$.
 
The natural morphism $\cH_{g,n+1}\ra\cH_{g,n}$ (forgetting the last marked point) is naturally
isomorphic to the universal $n$-punctured, genus $g$ curve over $\cH_{g,n}$ and the fiber over
any closed point $[C]\in\cH_{g,n}$ is diffeomorphic to $S_{g,n}$. Identifying its fundamental 
group with $\Pi_{g,n}$, we get, as above, a faithful topological monodromy representation:
$$\rho_{g,n}\co \pi_1(\cH_{g,n},[C])\ra\out(\Pi_{g,n}).$$
Instead, the faithfulness of the corresponding algebraic monodromy representation:
$$\hat{\rho}_{g,n}\co \hat{\pi}_1(\cH_{g,n},[C])\ra\out(\hP_{g,n}).$$
is a much deeper statement, equivalent to the congruence subgroup property for $H_{g,n}$.

The main result of this paper is that $\hat{\rho}_{g,n}$ is faithful for all $g$ and $n$ such that $g\geq 1$
and $n\geq 1$. 
In particular, we prove that the congruence subgroup property holds for the genus $2$ Teichm\"uller
modular group for $n\geq 1$ (the genus $0$ and $1$ cases have been proved by Asada in \cite{Asada}).

\section{The geometric profinite completion of $\GG_{g,n}$}\label{geometric pro}
Let us assume that the fundamental group $\Pi_{g,n}$ of $S_{g,n}$ has $P_{n+1}$ as base point.
For $2g-2+n>0$, the short exact sequence of topological 
fundamental groups, associated to the Serre fibration $\cM_{g,n+1}\ra\cM_{g,n}$, is then identified 
with the classical short exact sequence of modular groups
$$1\ra\Pi_{g,n}\ra\GG_{g,n+1}\ra\GG_{g,n}\ra 1,$$
while the corresponding short exact sequence of algebraic fundamental groups is identified with the 
short exact sequence
$$1\ra\hP_{g,n}\ra\hGG_{g,n+1}\ra\hGG_{g,n}\ra 1.$$
The action by inner automorphisms of $\hGG_{g,n+1}$ on its normal subgroup 
$\hP_{g,n}$ induces the representations $\tilde{\rho}_{g,n}\co\hGG_{g,n+1}\ra\aut(\hP_{g,n})$
and $\hat{\rho}_{g,n}\co\hGG_{g,n}\ra\out(\hP_{g,n})$. 

Let us mention here a fundamental result of Nikolov and Segal \cite{N-S} which asserts that any finite 
index subgroup of any topologically finitely generated profinite group $G$ is open. 
Since such a profinite group $G$ has also a basis of neighborhoods of the identity consisting of 
open characteristic subgroups, it follows that all automorphisms of $G$ are continuous and that
$\mbox{Aut}(G)$ is a profinite group as well. 
Let us then give the following definitions:

\begin{definition}\label{geometric}Let us define the profinite groups $\tGG_{g,n+1}$ 
and $\cGG_{g,n}$, for $2g-2+n>0$, to be, respectively, the image of $\tilde{\rho}_{g,n}$ in 
$\mbox{Aut}(\hat{\Pi}_{g,n})$ and of $\hat{\rho}_{g,n}$ in $\mbox{Out}(\hat{\Pi}_{g,n})$.
\end{definition}
 
By definition, there are natural maps with dense image $\GG_{g,n}\ra\cGG_{g,n}$ and 
$\GG_{g,n+1}\ra\tGG_{g,n+1}$, but it is a deep result by Grossman \cite{Grossman}
that these maps are also injective.

By Definition~\ref{geometric}, the representation $\tGG_{g,n+1}\ra\aut(\hP_{g,n})$, induced by the
action of inner automorphisms of $\tGG_{g,n+1}$ on its normal subgroup $\hP_{g,n}$, is injective.
Therefore, it holds:

\begin{proposition}\label{center 2} The center of $\tGG_{g,n+1}$ is trivial for $2g-2+n>0$.
\end{proposition}

Another consequence of Definition~\ref{geometric} is the following:

\begin{proposition}\label{shortexact}For $2g-2+n>0$, there is a natural short exact sequence:
$$1\ra\hP_{g,n}\ra\tGG_{g,n+1}\ra\cGG_{g,n}\ra 1.$$
In particular, $\cGG_{g,n}\equiv\hGG_{g,n}$ if and only if $\tGG_{g,n+1}\equiv\hGG_{g,n+1}$.
\end{proposition}
 
We then have the interesting corollary:

\begin{corollary}\label{center 1}If the congruence subgroup property holds for $\GG_{g,n}$, 
then $\hGG_{g,n+1}$ has trivial center.
\end{corollary}

A natural guess is that, for $2g-2+n>0$, the two profinite completions $\cGG_{g,n+1}$ and 
$\tGG_{g,n+1}$ of $\GG_{g,n+1}$ coincide. For $n>0$, this is a direct consequence of Theorem~2.2 in \cite{Matsu}:

\begin{theorem}[Matsumoto]\label{comparison}For $2g-2+n>0$ and $n\geq 1$, there is a natural isomorphism 
$\Phi\co\cGG_{g,n+1}\sr{\sim}{\ra}\tGG_{g,n+1}$. Hence, a short exact sequence: \,
$1\ra\hP_{g,n}\ra\cGG_{g,n+1}\ra\cGG_{g,n}\ra 1.$
\end{theorem}

The existence of a natural epimorphism $\cGG_{g,n+1}\ra\tGG_{g,n+1}$, for all $n\geq 0$, was already remarked in
the proof of Theorem~1 in \cite{Asada} and, as an immediate consequence, the genus $0$ case of the subgroup congruence property followed:

\begin{proposition}[Asada]\label{genus 0}For $n\geq 3$, it holds $\tGG_{0,n}\equiv\cGG_{0,n}\equiv\hGG_{0,n}$.
\end{proposition}
\begin{proof}The case $n=3$ is trivial, since $\GG_{0,3}=\{1\}$. The general case follows by 
Proposition~\ref{shortexact}, the epimorphism $\Phi\co\cGG_{0,n}\tura\tGG_{0,n}$ 
and induction on $n$.

\end{proof}

By Theorem~\ref{comparison}, for $2g-2+n>0$ and $n\geq 1$, we can define unambiguously, 
{\it the geometric profinite completion} of $\GG_{g,n}$ to be the group $\cGG_{g,n}$. 
An important corollary of the theorem is also the following:

\begin{corollary}\label{centerfree2}For $2g-2+n> 0$ and $n\geq 2$, the geometric profinite completion 
$\cGG_{g,n}$ has trivial center.
\end{corollary}

Asada in \cite{Asada} has proved the genus $1$ case of the subgroup congruence conjecture:

\begin{theorem}[Asada]\label{asada}It holds $\cGG_{1,n}\equiv\hGG_{1,n}$, for $n\geq 1$.
\end{theorem}

Let us remark, however, that the natural epimorphism $\hGG_{1,1}\ra\SL_2(\ZZ)$ is not injective
(see \S 8.8 in \cite{R-Z} for details). This is not a surprise, since $S_1$ is not a hyperbolic surface.

In the next section, in particular, we will also provide an alternative proof of Asada's Theorem.

\section{The hyperelliptic modular group}
\label{congruence}

In this section, we are going to prove the results announced in the introduction. The main feature of 
the moduli stack of $n$-pointed, genus $g$ smooth hyperelliptic complex curves $\cH_{g,n}$ is that
it can be described in terms of moduli of pointed genus $0$ curves. More precisely,
there is a natural $\Z/2$-gerbe $\cH_{g}\ra\cM_{0,[2g+2]}$, for $g\geq 2$, defined assigning, to a 
genus $g$ hyperelliptic curve $C$, the genus zero curve $C/\iota$, where $\iota$ is the hyperelliptic 
involution of $C$, labeled by the branch points of the cover $C\ra C/\iota$. In the genus $1$ case, 
there is a $\Z/2$-gerbe $\cM_{1,1}\ra\cM_{0,1[3]}$, where, by the notation "$1[3]$", we mean that 
one label is distinguished while the others are unordered. For $2g-2+n>0$,
there is also a natural representable morphism $\cH_{g,n+1}\ra\cH_{g,n}$, forgetting the 
$(n+1)$-th labeled point, which is isomorphic to the universal $n$-punctured curve over $\cH_{g,n}$.
So, the fiber above an arbitrary closed point $x\in\cH_{g,n}$ is diffeomorphic to $S_{g,n}$ and its
fundamental group is isomorphic to $\Pi_{g,n}$. These morphisms induce, on 
topological fundamental groups, the short exact sequences, for $g\geq 2$:
$$
1\ra\Z/2\ra H_g\ra\GG_{0,[2g+2]}\ra 1\,\,\,\,
\mbox{ and }\,\,\,\,
1\ra\Pi_{g,n}\ra H_{g,n+1}\ra H_{g,n}\ra 1.
$$
Similarly, for the algebraic fundamental groups, there are short exact sequences:
$$
1\ra\Z/2\ra \hH_g\ra\hGG_{0,[2g+2]}\ra 1\,\,\,\,
\mbox{ and }\,\,\,\,
1\ra\hP_{g,n}\ra \hH_{g,n+1}\ra \hH_{g,n}\ra 1.
$$

The outer representation $\hat{\rho}_{g,n}\co\hH_{g,n}\ra\out(\hP_{g,n})$, induced by the last of the
above short exact sequences, is the algebraic monodromy representation of the punctured universal
curve over $\cH_{g,n}$. As already remarked, the congruence subgroup property for $H_{g,n}$ is
equivalent to the faithfullness of $\hat{\rho}_{g,n}$.

Let us prove some general properties of the groups $H_{g,n}$.
For definitions and elementary properties of good groups, we refer
to exercise 1 in Section~2.6 of \cite{Serre}. From the above exact sequences, it then follows immediately:

\begin{proposition}\label{good}For $2g-2+n>0$ and $g\geq 1$, the group $H_{g,n}$ is good.
\end{proposition}

It is well known that the centralizer of a finite index subgroup $U$ of $H_{g,n}$, 
for $g\geq 2$ and $n=0$ or $g=1$ and $n=1$, is spanned by the hyperelliptic involution $\iota$ while
it is trivial for $g\geq 2$ and $n\geq 1$ or $g=1$ and $n\geq 2$. 
An analogous statement holds for the profinite completion $\hH_{g,n}$.

\begin{proposition}\label{centerfree}Let $U$ be an open subgroup of $\hH_{g,n}$, for $2g-2+n>0$.
Then, for $g\geq 2$ and $n=0$ or $g=1$ and $n=1$, the centralizer of $U$ in $\hH_{g,n}$ is 
spanned by the hyperelliptic involution. In all the other cases, the centralizer of $U$ in $\hH_{g,n}$ 
is trivial.
\end{proposition}

\begin{proof}Let us consider first the cases $g\geq 2$ and $n=0$ or $g=1$ and $n=1$.
It is clearly enough to prove that for any open subgroup $U$ of $\hH_{g,n}$, which
contains the hyperelliptic involution $\iota$, the center $Z(U)$ is equal to the subgroup 
spanned by $\iota$.

The center of any open subgroup of $\hGG_{0,[2g+2]}$ is trivial. From the exact sequences:
$$1\ra\Z/2\cdot\iota\ra\hH_g\ra\hGG_{0,[2g+2]}\ra 1\,\,\,\,\mbox{ and }\,\,\,\,
1\ra\Z/2\cdot\iota\ra\hGG_{1,1}\ra\hGG_{0,[4]},$$
it then follows that $Z(U)=\langle\iota\rangle$. 

For the cases $g\geq 2$ and $n\geq 1$ or $g=1$ and $n\geq 2$, we have to prove
that the center is trivial for any open subgroup $U$ of $\hH_{g,n}$. 
By induction on $n$, thanks to the short exact sequences:
$$1\ra\hP_{g,n-1}\ra\hH_{g,n}\ra \hH_{g,n-1}\ra 1,$$
it is enough to prove the proposition for the cases $g\geq 2$, $n=1$ and $g=1$, $n=2$. 

From the above short exact sequence, we then see that the center $Z(U)$, if non-trivial, projects 
to the subgroup of $\hH_{g,n-1}$ spanned by the hyperelliptic involution. 

In this case, the subgroup $Z(U)\cdot\hP_{g,n-1}$ of $\hH_{g,n}$ would be generated by a hyperelliptic 
involution $\mu$ in $H_{g,n}$ and $\hP_{g,n-1}$. So, $Z(U)$ would be generated by a conjugate
$f\mu f^{-1}$ for some $f\in\hP_{g,n-1}$. Let $U':=fUf^{-1}$, then it is clear that 
$Z(U')=\langle\mu\rangle$.  

Hence, such $\mu$ would commute with the elements of the finite 
index subgroup $U'\cap\Pi_{g,n-1}$ of $\Pi_{g,n-1}$. By the simple topological description of
a hyperelliptic involution, there is a simple loop $\gm\in\Pi_{g,n-1}$ such that $\mu(\gm)=\gm^{-1}$. 

As already remarked in \S~\ref{geometric pro}, if we identify $\hP_{g,n-1}$ with its image in
$\hH_{g,n}$, then it holds $\mu(\gm)=\mu\gm\mu^{-1}$.
For some $k>0$, it also holds $\gm^k\in U'\cap\Pi_{g,n-1}$ and then:
$$\gm^{-k}=\mu(\gm^k)=\mu\gm^k\mu^{-1},$$
which contradicts the fact that $Z(U')=\langle\mu\rangle$. Therefore, it holds $Z(U)=\{1\}$.

\end{proof}

We call a finite index subgroup $H^\ld$ of $H_{g,n}$ a {\it level} of $H_{g,n}$ and the corresponding 
\'etale cover $\cH^\ld\ra\cH_{g,n}$ {\it a level structure} over $\cH_{g,n}$.
Geometric levels of $H_{g,n}$ are defined by means of the monodromy 
representation $\rho\co H_{g,n}\ra\out(\Pi_{g,n})$. For a characteristic subgroup $\Pi^\ld$ of $\Pi_{g,n}$, 
the geometric level $H^\ld$ is defined to be the kernel of the induced representation 
$\rho_\ld\co H_{g,n}\ra\out(\Pi_{g,n}/\Pi^\ld)$. The abelian level 
$H(m)$ of order $m\geq 2$ is then defined to be the kernel of the representation 
$\rho_{(m)}\co H_{g,n}\ra\Sp_{2g}(\Z/m)$ and we let 
$\cH^{(m)}$ be the corresponding abelian level structure. 

There is a standard procedure to simplify the structure of an algebraic stack $X$ by 
erasing a generic group of automorphisms $G$ (see, for instance, \cite{Romagny}). The algebraic 
stack thus obtained is usually denoted by $X\!\!\fatslash G$. So, the natural map 
$\cH_g\ra\cM_{0,[2g+2]}$ yields an isomorphism 
$\cH_g\!\!\fatslash\langle\iota\rangle\cong\cM_{0,[2g+2]}$. A natural question is then which level
structure over $\cH_g$ corresponds to the Galois \'etale cover $\cM_{0,2g+2}\ra\cM_{0,[2g+2]}$.

\begin{proposition}\label{ordered}For $g\geq 2$, there is a natural isomorphism 
$\cH^{(2)}\!\!\!\fatslash\langle\iota\rangle\cong\cM_{0,2g+2}$.
\end{proposition}

\begin{proof}The groups $H_g/\langle\iota\rangle$ and $\GG_{0,[2g+2]}$ 
are naturally isomorphic. By means of this 
isomorphism, the normal subgroup $\GG_{0,2g+2}\lhd\GG_{0,[2g+2]}$ identifies with the subgroup 
of $H_g/\langle\iota\rangle$ spanned by squares of Dehn twists along non-separating s.c.c. on 
$S_{g,n}$. Squares of Dehn twists, all act trivially on homology with $\Z/2$-coefficients. Therefore,
$\GG_{0,2g+2}$ identifies with a normal finite index subgroup of $H(2)/\langle\iota\rangle$.
So, there are a natural \'etale morphism $\Phi\co\cM_{0,2g+2}\ra\cH^{(2)}\!\!\!\fatslash\langle\iota\rangle$
and a commutative diagram with exact rows:
$$\begin{array}{ccccccl}
1\ra&\GG_{0,2g+2}&\ra&\GG_{0,[2g+2]}&\ra&\Sigma_{2g+2}&\ra 1\\
&\cap&&\|\wr&&\da{\sst \rho}&\\
1\ra&H(2)/\langle\iota\rangle&\ra&H_g/\langle\iota\rangle&\ra&\operatorname{PGL}_{2g}(\Z/2).&
\end{array}$$
At this point, observe that the representation $\rho\co\Sigma_{2g+2}\ra\operatorname{PGL}_{2g}(\Z/2)$ is 
induced by the permutation of $2g+2$ points in general position in the projective space 
$\PP_{\Z/2}^{2g-1}$ and so is faithful. Thus, the injection $\GG_{0,2g+2}\hookra H(2)/\langle\iota\rangle$ 
is actually an isomorphism and then $\Phi$ is an isomorphism as well.

\end{proof}

\begin{remark}Likewise, it is not hard to prove that, for the abelian level structure 
$\cM^{(2)}$ over $\cM_{1,1}$, there is a natural isomorphism 
$\cM^{(2)}\!\!\!\fatslash\langle\iota\rangle\cong\cM_{0,4}$,
where $\iota$ here denotes the generic elliptic involution.

From now on, we will mostly stick to moduli spaces of hyperelliptic curves $\cH_{g,n}$, with $g\geq 2$, 
and leave to the reader the formulation and the proof of the analogous statements for $g=1$, $n\geq 1$.
\end{remark}

Let $\cC_g\ra\cH_g$, for $g\geq 2$, be the universal curve. Removing Weierstrass points from its fibers,
we obtain a $(2g+2)$-punctured, genus $g$ curve $\cC_0\ra\cH_{g}$. A weak 
version of the congruence subgroup property for $H_g$ is then the assertion that the algebraic 
monodromy representation, associated to $\cC_0\ra\cH_{g}$, is faithful:
$$\hat{\rho}_0\co\hat{\pi}_1(\cH_g,x)\hookra\out(\hat{\pi}_1(C_0,\ol{x})),$$ 
where $C_0$ is the fiber of $\cC_0\ra\cH_g$ over the closed point $x$. Let us show how 
this assertion reduces to Corollary~\ref{genus 0}.

Let us denote by $\cC^\ld\ra\cH^\ld$ the pull-back of the universal curve $\cC_g\ra\cH_g$ to
the level structure $\cH^\ld\ra\cH_g$ and by $\cC_0^\ld\ra\cH^\ld$ the pull-back of the 
punctured curve $\cC_0\ra\cH_g$. 

By Proposition~\ref{ordered}, there is a natural \'etale Galois morphism 
$\cH^{(4)}\ra\cM_{0,2g+2}$ which is also representable, since $\iota\notin H(4)$. 
Let ${\cal R}\ra\cH^{(4)}$ be the pull-back of the universal
$(2g+2)$-punctured, genus $0$ curve $\cM_{0,2g+3}\ra\cM_{0,2g+2}$.
There is then a commutative diagram:
$$\begin{array}{ccc}
\cC_0^{(4)}&\sr{\psi}{\ra}&{\cal R}\,\,\,\\
&\searrow&\da\\
&& \cH^{(4)},
\end{array}
$$
where $\psi$ is the \'etale, degree $2$ map which, fiberwise, is the quotient by the hyperelliptic 
involution. The algebraic monodromy representation
$\hat{\pi}_1( \cH^{(4)},a)\ra\out(\hat{\pi}_1({\cal R}_a,\ol{a}))$, associated to the rational curve
${\cal R}\ra\cH^{(4)}$, is faithful by Corollary~\ref{genus 0}. 
Then, by Lemma~8 in \cite{Asada}, the algebraic monodromy representation 
$\hat{\pi}_1( \cH^{(4)},a)\ra\out(\hat{\pi}_1(C_0^{(4)},\tilde{a}))$, associated to the curve 
$\cC_0^{(4)}\ra\cH^{(4)}$, is faithful as well, where $C_0$ denotes the fiber over the 
closed point $a$. This immediately implies the faithfulness of the representation $\hat{\rho}_0$.

We can now state and prove the main result of the paper:

\begin{theorem}\label{congruencepro} Let $\cH_{g,n}$, for $2g-2+n>0$ and $g\geq 1$, be the moduli
stack of $n$-pointed, genus $g$ hyperelliptic complex curves. For $n\geq 1$,   
the universal algebraic monodromy representation 
$\hat{\rho}_{g,n}\co\hat{\pi}_1( \cH_{g,n})\ra\out(\hP_{g,n})$,
associated to the universal $n$-punctured, genus $g$ 
hyperelliptic curve $\cH_{g,n+1}\ra\cH_{g,n}$, is faithful.
\end{theorem}

\begin{proof}The proof of Theorem~\ref{congruencepro} consists of two steps. In the first, we show 
that the faithfulness of $\hat{\rho}_{g,n}$, for a given $g\geq 1$ and all $n\geq 1$, 
can be deduced from that of $\hat{\rho}_{g,n'}$, for any given $n'\geq 1$. 
In the second, we prove that $\hat{\rho}_{g,2g+2}$ is faithful for all $g\geq 1$. 
The first step is accomplished, by induction, in the following lemma: 

\begin{lemma}\label{n=n+1}Let $g\geq 1$ and $n\geq 1$. Then, the monodromy representation 
$\hat{\rho}_{g,n}$ is faithful if and only if $\hat{\rho}_{g,n+1}$ is.
\end{lemma}
\begin{proof}By Theorem~\ref{comparison}, for $n\geq 1$, there is a commutative diagram with exact rows:
$$\begin{array}{ccccccc}
1\ra&\hP_{g,n}&\ra&\hH_{g,n+1}&\ra&\hH_{g,n}&\ra 1\\
&\|\,\,\,\,\,&&\da{\st \hat{\rho}_{g,n+1}}\,&&\,\,\da{\st \hat{\rho}_{g,n}}&\\
1\ra&\hP_{g,n}&\ra&\cGG_{g,n+1}&\ra&\cGG_{g,n}&\ra 1
\end{array}$$
and the lemma follows immediately.

\end{proof}

\begin{lemma}\label{final step}For $g\geq 1$, the algebraic monodromy representation
$\hat{\rho}_{g,2g+2}$ is faithful.
\end{lemma}

\begin{proof}Here, as usual, for notational reason, we assume $g\geq 2$ and leave to the reader 
the transposition of the argument to the genus $1$ case.

The universal curve $\cC^{(2)}\ra\cH_g^{(2)}$ is endowed with $2g+2$ ordered sections, 
corresponding to the Weierstrass points on the fibers.   
So, by the universal property of $\cH_{g,n}^{(2)}$, there is a morphism 
$s\co\cH_g^{(2)}\ra\cH_{g,2g+2}^{(2)}$ which is a section of the natural projection 
$p\co\cH_{g,2g+2}^{(2)}\ra\cH_g^{(2)}$ (forgetting the labels). The morphism $p$ is 
smooth and its fiber above a closed point $[C]\in\cH_g$ is the configuration space 
of $2g+2$ points on the curve $C$. Let us denote by $S_g(n)$ the configuration space of
$n$ points on the compact Riemann surface $S_g$ and by $\Pi_g(n)$ its fundamental group.
Then, all fibers of $p$ above closed points of $\cH_g$ are diffeomorphic to $S_g(n)$. Therefore,
the fundamental group $H_{g,2g+2}(2)$ of $\cH_{g,2g+2}^{(2)}$ fits in the short exact sequence:
$$1\ra\Pi_g(2g+2)\ra H_{g,2g+2}(2)\ra H_g(2)\ra 1,$$
which is split by $s_\ast\co H_g(2)\ra H_{g,2g+2}(2)$. Moreover, since the space $S_g(n+1)$ is 
fibered in $n$-punctured, genus $g$ curves over $S_g(n)$, for all $n\geq 0$, there is
a short exact sequence:
$$1\ra\Pi_{g,n}\ra\Pi_g(n+1)\ra\Pi_g(n)\ra 1.$$

From Theorem~\ref{comparison} and a simple induction on $n$, it follows that the profinite 
completion $\hP_g(n)$ embeds in $\cGG_{g,n}$ (this is essentially the same argument of Asada in 
Theorem~1, \cite{Asada}, where this was first proved). Therefore, passing to profinite completions, 
we get the short exact sequences:
$$\begin{array}{c}
1\ra\hP_g(2g+2)\ra \hH_{g,2g+2}(2)\ra \hH_g(2)\ra 1,\\
\\
1\ra\hP_{g,n}\ra\hP_g(n+1)\ra\hP_g(n)\ra 1.
\end{array}$$
The former is split by $\hat{s}_\ast\co \hH_g(2)\ra \hH_{g,2g+2}(2)$. So there is an isomorphism:
$$\hH_{g,2g+2}(2)\cong\hP_g(2g+2)\rtimes\hH_g(2).$$

In order to prove that the algebraic monodromy representation 
$\hat{\rho}_{g,2g+2}$ is faithful, it is enough to show that this holds for its restriction to
$\hH_{g,2g+2}(2)$, which we denote also by $\hat{\rho}_{g,2g+2}$. But we have already seen that 
$\hat{\rho}_{g,2g+2}\circ\hat{s}_\ast=\hat{\rho}_0\co \hH_g(2)\ra\out(\hP_{g,2g+2})$ is faithful and, 
as remarked above, the restriction of $\hat{\rho}_{g,2g+2}$ to the normal subgroup 
$\hP_g(2g+2)$ of $\hH_{g,2g+2}(2)$ is faithful as well. So, Lemma~\ref{final step} follows, 
if we prove that:
$$\hat{\rho}_{g,2g+2}(\hP_g(2g+2))\cap\hat{\rho}_{g,2g+2}(\hat{s}_\ast( \hH_g(2)))=\{1\}
\,\,\,\,\,\,\,\,\,(\ast).$$

The subgroup $s_\ast(H_g(2))$ of $H_{g,2g+2}(2)$ centralizes the 
hyperelliptic involution $s_\ast(\iota)\in H_{g,2g+2}(2)$. Passing to profinite completions, 
the subgroup $\hat{s}_\ast(\hH_g(2))$ of $\hH_{g,2g+2}(2)$ then centralizes the 
hyperelliptic involution $s_\ast(\iota)\in \hH_{g,2g+2}(2)$. 
It is clear that $\hat{\rho}_{g,2g+2}(\hat{s}_\ast(\iota))\neq 1$. Hence, since
$\hP_g(2g+2)$ is torsion free:
$$\hat{\rho}_{g,2g+2}(\hP_g(2g+2)\cdot s_\ast(\iota))\cong\hP_g(2g+2)\cdot s_\ast(\iota)
\cong\hP_g(2g+2)\rtimes\Z/2.$$
All elements of $\hat{\rho}_{g,2g+2}(\hat{s}_\ast( \hH_g(2)))$ commute with 
$\hat{\rho}_{g,2g+2}(s_\ast(\iota))$. So, in order to prove the identity $(\ast)$, 
it is enough to show that no element of $\hat{\rho}_{g,2g+2}(\hP_g(2g+2))$ does.

From item $(ii)$ of Lemma~2.1 in \cite{Mochi}, it follows that a primitive finite subgroup of
the algebraic fundamental group of a hyperbolic orbi-curve is self-normalizing.

For all $0\leq n\leq 2g+1$, a given hyperelliptic involution $\iota'\in H_{g,n+1}$ and $\hP_{g,n}$
span inside of $\hH_{g,n+1}$ a group isomorphic to the algebraic fundamental group
of an $n$-punctured, genus $g$ hyperelliptic orbi-curve $[C/\iota']$. In particular, by 
Lemma~2.1 in \cite{Mochi}, there is no element of $\hP_{g,n}$ with which $\iota'$ commutes.

The short exact sequences $1\ra\hP_{g,n}\ra\hP_g(n+1)\ra\hP_g(n)\ra 1$ and
a simple induction on $n\geq 0$ then imply that $s_\ast(\iota)$ does not commute with any 
given element of $\hP_g(2g+2)$, as claimed above. This completes the proof of 
Lemma~\ref{final step} and then that of Theorem~\ref{congruencepro}.

\end{proof}
\end{proof}

\subsection*{Acknowledgements}
This work was basically conceived during my stay in RIMS in the Fall of 2007. It was a wonderful 
and enriching experience to be there. I had the unique opportunity to speak with the world best experts
in the area of anabelian geometry, to which I was trying to contribute in that period and
in which, I think, the present paper should be collocated.

Accordingly, I thank Makoto Matsumoto, who made available his personal research funds 
to finance my stay there, and Shinichi Mochizuki, who invited me to RIMS.
I also thank them and Akio Tamagawa, for the warm hospitality and the many 
mathematical discussions we had. Finally, I thank the referee for the many useful comments 
he made on preliminary versions of this work.

\bigskip

\noindent Address:\, Escuela de Matem\'atica, Universidad de Costa Rica, San Jos\'e, Costa Rica.
\\
E--mail:\, marco.boggi@gmail.com

\end{document}